\documentclass[a4paper,10pt]{article}

\usepackage[english]{babel}
\usepackage{german, umlaute, amsmath, amssymb , latexsym, theorem,
  epic, epsfig, rotating}

\newcommand{\xn}{\ensuremath{x_{n}}}

\newcommand{\nk}{\ensuremath{n_{k}}}

\newcommand{\xnull}{\ensuremath{x_{0}}}

\newcommand{\neins}{\ensuremath{n_{1}}}

\newcommand{\xeins}{\ensuremath{x_{1}}}

\newcommand{\ntel}{\ensuremath{\frac{1}{n}}}

\newcommand{\qtel}{\ensuremath{\frac{1}{q}}}

\newcommand{\halb}{\ensuremath{\frac{1}{2}}}
\newcommand{\drittel}{\ensuremath{\frac{1}{3}}}
\newcommand{\viertel}{\ensuremath{\frac{1}{4}}}

\newcommand{\dreiviertel}{\ensuremath{\frac{3}{4}}}

\newcommand{\alphlist}{\begin{list}{(\alph{enumi})}{\usecounter{enumi}}}
\newcommand{\romanlist}{\begin{list}{(\roman{enumi})}{\usecounter{enumi}}}
\newcommand{\listend}{\end{list}}

\newcommand{\ldot}{\ensuremath{\textbf{.}}}
\newcommand{\ld}{\ensuremath{,\ldots,}}

\newcommand{\ra}{\ensuremath{\rightarrow}}
\newcommand{\equi}{\ensuremath{\Leftrightarrow}}
\newcommand{\follows}{\ensuremath{\Rightarrow}}

\newcommand{\N}{\ensuremath{\mathbb{N}}} 
\newcommand{\R}{\ensuremath{\mathbb{R}}}
\newcommand{\Z}{\ensuremath{\mathbb{Z}}}

\newcommand{\kreis}{\ensuremath{\mathbb{T}^{1}}}
\newcommand{\ntorus}[1][2]{\ensuremath{\mathbb{T}^{#1}}}

\newcommand{\nLim}{\ensuremath{\lim_{n\rightarrow\infty}}}

\newcommand{\jLim}{\ensuremath{\lim_{j\rightarrow\infty}}}

\newcommand{\nfolge}[1]{\ensuremath{(#1)_{n\in\mathbb{N}}}}

\newcommand{\qed}{{\raggedleft $\Box$ \\}}
\newcommand{\proof}{\textit{Proof:}}

\newcommand{\graph}{\ensuremath{\mathrm{graph}}}

\newcommand{\hin}{\ensuremath{h^{-1}}}

\newcommand{\pie}{\ensuremath{\pi_{1}}}
\newcommand{\piz}{\ensuremath{\pi_{2}}}

\newcommand{\piin}{\ensuremath{\pi^{-1}}}

\newcommand{\ncup}{\ensuremath{\bigcup_{n\in\N}}}

\newcommand{\ikcap}{\ensuremath{\bigcap_{i=1}^k}}

\newcommand{\twoscriptsum}[2]{ \ensuremath{\sum_{\begin{array}{c}
        \scriptstyle #1 \\ \scriptstyle #2 \end{array} }}}

\newcommand{\insum}{\ensuremath{\sum_{i=1}^n}}

\newcommand{\inergsum}{\ensuremath{\sum_{i=0}^{n-1}}}

\newcommand{\Thom}{\ensuremath{{\cal T}_{\textrm{hom}}}}
\newcommand{\Tdiff}{\ensuremath{{\cal T}_{\textrm{diff}}}}
\newcommand{\Tbv}{\ensuremath{{\cal T}_{\textrm{BV}}}}

\newcommand{\muth}{\ensuremath{\mu_{\theta}}}
\newcommand{\muthom}{\ensuremath{\mu_{\theta+\omega}}}

\newcommand{\Ath}{\ensuremath{A_{\theta}}}

\newcommand{\thx}{\ensuremath{(\theta,x)}}
\newcommand{\thom}{\ensuremath{\theta + \omega}}

\newcommand{\omtil}{\ensuremath{\tilde{\omega}}}

\newcommand{\phith}{\ensuremath{\varphi(\theta)}}
\newcommand{\pqphith}{\ensuremath{\displaystyle (\varphi^i_j(\theta))^{
 \scriptscriptstyle 1\leq i \leq p}_{ \scriptscriptstyle 1 \leq j \leq
q }  } }

\newcommand{\psith}{\ensuremath{\psi(\theta)}}
 
\newcommand{\pqphiset}{\ensuremath{\{
\varphi_1^1(\theta), \ldots , \varphi^p_q(\theta) \} }}

\newcommand{\muphi}{\ensuremath{\mu_{\varphi}}}

\newcommand{\gamhat}{\ensuremath{\widehat{\gamma}}}

\newcommand{\Tth}{\ensuremath{T_{\theta}}}
\newcommand{\Tnth}{\ensuremath{T_{\theta}^{n}}}

\newcommand{\Tthx}{\ensuremath{T_{\theta}(x)}}
\newcommand{\Tnthx}{\ensuremath{T_{\theta}^{n}(x)}}

\newcommand{\That}{\ensuremath{\widehat{T}}}
\newcommand{\Thatth}{\ensuremath{\widehat{T}_\theta}}
\newcommand{\Tbar}{\ensuremath{\overline{T}}}
\newcommand{\thhat}{\ensuremath{\widehat{\theta}}}
\newcommand{\xhat}{\ensuremath{\widehat{x}}}

\newcommand{\phihat}{\ensuremath{\widehat{\varphi}}}

\theoremstyle{break}
\newtheorem{definition}{Definition}[section]
\newtheorem{thm}[definition]{Theorem}
\newtheorem{thmdef}[definition]{Theorem and Definition}
\newtheorem{lem}[definition]{Lemma}
  
\newtheorem{prop}[definition]{Proposition}

\theorembodyfont{\rmfamily}
\newtheorem{bem}[definition]{Remark}    
\newtheorem{example}[definition]{Example}

\numberwithin{equation}{section}

\title{The Denjoy type-of argument for quasiperiodically forced circle
  diffeomorphisms} \author{Tobias H. J\"ager and Gerhard Keller \\
  \small Mathematisches
  Institut,  \\
  \small Friedrich-Alexander-Universit\"at Erlangen-N\"urnberg,
  Germany \\ \small MSC Classification Numbers: 37C70, 37C60, 37H15 \\
}

\begin{document}
\setlength{\topmargin}{-1cm}
\setlength{\oddsidemargin}{-0.03\textwidth}  

\maketitle
%
%
%                               Abstract
%
\begin{abstract}
  We carry the argument used in the proof of the Theorem of Denjoy
  over to the quasiperiodically forced case. Thus we derive that if a
  system of quasiperiodically forced circle diffeomorphisms with
  bounded variation of the derivative has no invariant graphs with a
  certain kind of topological regularity, then the system is
  topologically transitive.
\end{abstract}

\section{Introduction}  \label{Intro}

When studying quasiperiodically forced circe maps, it is an obvious
question to ask to what extent the results for unperturbed circle
homeomorphisms and diffeomorphisms carry over.  Invariant graphs serve
as a natural analogue for fixed or periodic points. Unfortunately, a
system without invariant graphs is not necessarily conjugated to an
irrational torus translation (at least not if the conjugacy is
required to respect the skew product structure; a simple
counterexample will be given in Section \ref{Furstenberg}).  However,
there is a result by Furstenberg (\cite{furstenberg:1961}), which may
be considered as an anologue to the Poincar\'{e} classification for
circle homeomorphisms in a measure-theoretic sense: It states that for
a system of quasiperiodically forced circle maps either there exists
an invariant graph, then all invariant ergodic measures are associated
to invariant graphs, or the system is uniquely ergodic with respect to
an invariant measure which is continuous on the fibres. In the latter
case the system will be isomorphic to a skew translation of the torus
by an isomorphism with some additional nice properties. We will
briefly discuss Furstenbergs results after we have introduced the
concept of invariant graphs in Section \ref{Furstenberg}.

The (fibrewise) rotation number can be defined in the same way as the
rotation number of a circle homeomorphism, but Herman already gave
examples where invariant graphs exist in combination with arbitrary
rotation numbers (see \cite{herman:1983}). Thus the connection between
the existence of fixed or periodic points and rational rotation
numbers seems to break down in the forced case. However, this is
different if we require a certain amount of topological regularity of
the invariant graphs, in particular when the invariant graphs are
continuous. Then the structure of an invariant graph already determines the
fibrewise rotation number, and this is described in Section \ref{Rotnum}.
It might be considered well-known folklore, but as we knew of no
apropriate reference and need a slight generalization afterwards it
seemed apropriate to include some details.

Section \ref{Denjoy} contains our main result, namely Theorem
\ref{thm:denjoy}. $\Tbv$ denotes the class of systems where the
variation of the logarithm of the derivative of the fibre maps is
integrable. The concept of a regular $p,q$-invariant graph and its
implications for the rotation numbers are introduced in Sections
\ref{Furstenberg} and \ref{Rotnum}, as mentioned. The statement we
derive then reads as follows: :
\begin{quote} \textit{
    If $T \in \Tbv$ is not topologically transitive, then there exists
    a regular $p,q$-invariant graph for $T$. In particular, $\rho_T$
    depends rationally on $\omega$, i.e.\ $\rho_T = \frac{k}{q}\omega
    + \frac{l}{pq} \bmod 1$ for suitable integers $k$ and $l$. }
\end{quote}
The question whether such systems are minimal (as 
in the unforced case) has to be left open here. However, if there
exists a minimal strict subset of the circle then some further
conclusions about its structure can be drawn.

%
%
%                             Invariant graphs, Furstenberg
%
%

\section{Invariant graphs and invariant measures} \label{Furstenberg}

We study quasiperiodically forced circle homeomorphisms and
diffeomorphisms, i.e.\ continuous maps of the form
\begin{equation}                               
  \label{generalsystem}
  T : \ntorus \ra \ntorus \ , \ \thx \mapsto (\thom,\Tthx) \ ,
\end{equation}
where the \textit{fibre maps} \Tth\ are either orientation-preserving circle
homeomorphisms or orientation-preserving circle diffeomorphisms with
the derivative $D\Tth$ depending continuously on \thx. To ensure all
required lifting properties we additionally assume that $T$ is
homotopic to the identity on \ntorus.  The classes of
such systems will be denoted by \Thom\ and \Tdiff\ respectively.  In
all of the following, $m$ will be the Lebesgue-measure on \kreis,
$\lambda$ the Lebesgue-measure on \ntorus\ and
$\pi_{i} : \ntorus \ra \kreis \ (i=1,2)$ the projection to the
respective coordinate. When considering fibre maps of iterates of $T$
or their inverses, we use the convention
$\Tnth := (T^n)_{\theta}\ \forall n\in\Z$. Finally, by $\ntorus[n]_*$
we denote the set of all points $(x_1 \ld x_n)\in \ntorus[n]$ with
$x_i \neq x_j \ \forall i,j: 1 \leq i < j \leq n$.
 
\ \\
\textit{Invariant graphs.} Due to the aperiodicity of the forcing
rotation, there cannot be any fixed or periodic points for a system of
the form (\ref{generalsystem}). The simplest invariant objects will
therefore be invariant graphs. In contrast to quasiperiodically forced
monotone interval maps, where such invariant graphs are always
one-valued functions, these might be multi-valued when circle maps are
considered and the definition therefore needs a little bit more care.
We will usually not distinguish between invariant graphs as functions
and as point sets. This might seem a bit confusing at some times, but
it is very convenient at others.
\begin{definition}[$p,q$-invariant graphs]
  \label{def:invgraph}
  Let $T \in \Thom$, $p,q \in \N$. A $p,q$-invariant graph is a
  measurable function $\varphi : \kreis \ra \ntorus[pq]_*\ , \ \theta
  \mapsto \pqphith$ with
  \begin{equation}
    \label{eq:invgraph}
    \Tth(\pqphiset ) \ = 
    \ \{ \varphi_1^1(\thom) , \ldots , \varphi^p_q(\thom) \} 
  \end{equation}
  for $m$-a.e. $\theta \in \kreis$, which satisfies the following:
  \romanlist
  \item $\varphi$ cannot be decomposed (in a measurable way) into
  disjoint subgraphs $\varphi^1 \ld \varphi^m \ (m\in\N)$ which
  also satisfy (\ref{eq:invgraph}).
  \item $\varphi$ can be decomposed into $p$ $p$-periodic $q$-valued
    subgraphs $\varphi^1,\ldots,\varphi^p$.
  \item The subgraphs $\varphi^1 \ld \varphi^p$ cannot further be
    decomposed into invariant or periodic subgraphs.
  \listend
  If $p=q=1$, $\varphi$ is called a simple invariant graph.  The point
  set $\Phi := \{ (\theta,\varphi^i_j(\theta)) \mid \theta \in \kreis,\
  1\leq i \leq p,\ 1 \leq j \leq q \}$ will also be called an
  invariant graph, but labeled with the corresponding capital letter.
  An invariant graph is called continuous if it is continuous as a
  function $\kreis \ra \ntorus[pq]_*$.
\end{definition}
Note that by this convention an invariant graph
is a minimal object in the sense that it is not decomposable into
smaller invariant parts. Therefore, the union of two or more invariant
graphs will not be called an invariant graph again. On the other hand,
it is always possible to decompose an invariant set which is the graph
of a $n$-valued function into the disjoint union of invariant graphs:

\begin{lem}
  \label{lem:decomposition}
Let $T\in\Thom$ and suppose $F : \kreis \ra \ntorus[n]$ satisfies 
\[\Tth(\{
F_1(\theta), \ldots , F_n(\theta) \})\  = \ \{
F_1(\thom),\ldots,F_n(\thom) \} \ \ m\textrm{-a.s.\ .}
\]
Then $\textrm{graph}(F)$ can be
decomposed in exactly one way (modulo permutation) into $k$ disjoint
$p,q$-invariant graphs $\Phi^i$, i.e.\ $\textrm{graph}(F) =
\bigcup_{i=1}^k \Phi^k$, with $kpq = n$.
\footnote{\label{foot:randomset} For the sake of completeness, we
  should also mention the following: Suppose a measurable set consists
  of $n$ points on each fibre. One might ask whether such a set is
  always the graph of a $n$-valued measurable function.  Prop.\ 1.6.3
  in \cite{arnold:1998} provides a positive answer to this, namely
  \alphlist
\item
  If $A\subseteq \ntorus$ is a measurable set with $\# \Ath = n \
  \forall \theta \in \kreis$, then $A$ is the graph of a $n$-valued
  measurable function $F$. 
\item Let $A \subseteq \ntorus$ be a measurable set which is
  $T$-invariant and satisfies $\# \Ath < \infty$ $m$-a.s.~. Then $\#
  \Ath$ is $m$-a.s.\ constant and there exists a multi-valued function
  $F$ such that $\Ath = (\graph(F))_\theta$ $m$-a.s.~.  
\listend}
\end{lem}
\proof \\
\underline{Uniqueness:} Any two invariant graphs are either disjoint
or equal on $m$-a.e.\ fibre, otherwise their intersection would define
an invariant subgraph. Thus if $\Psi^1,\ldots,\Psi^l$ is another
decomposition of $\textrm{graph}(F)$ into invariant graphs, then every
$\Psi^j$ must be equal to some $\Phi^i$ (as it cannot be disjoint to
all of them). This immediately implies the uniqueness of the
decomposition.
\ \\
\underline{Existence:} $\textrm{graph}(F)$ is either an invariant
graph itself or can be decomposed into subgraphs which are invariant
as point sets. The same is true for
any such subgraph and after at most $n$ steps this yields an invariant
graph $\hat{\varphi}$. This is either a $1,q$-periodic invariant
graph, or it contains a periodic subgraph. Again, after a finite
number of steps this yields some $p$-periodic graph $\varphi^1 =
\varphi_1^1,\ldots,\varphi_q^1$ which is not further decomposable into
invariant or periodic subgraphs. W.l.o.g.\ we can assume that the
points $\varphi^1_j$ are ordered in \kreis, i.e.\ $\varphi^1_1(\theta) <
\varphi^1_2(\theta) < \ldots < \varphi^1_q(\theta) <
\varphi^1_1(\theta)$. Now for $m$-a.e.\ $\theta\in\kreis$
all the intervals $[\varphi^1_i(\theta),\varphi^1_{i+1}(\theta)]$
contain the same number of points $F_j(\theta)$, otherwise it would be
possible to define an $T^p$-invariant subgraph of $\varphi^1$. Thus,
by setting $\varphi_i^l(\theta) := $ $l$th point of
$\{F_1(\theta),\ldots,F_n(\theta) \}$ in
$[\varphi^1_i(\theta),\varphi^1_{i+1}(\theta)]$, the required
decomposition of $\textrm{graph}(F)$ into $p,q$-invariant graphs can
be defined.

\qed

\ \\
Of course, an example of a $p,1$-invariant graph can always be given
by choosing a circle homeomorphism $f$ with a periodic point of period
$p$ and taking $\Tth = f \ \forall \theta \in \kreis$. Then the
periodic orbit of $f$ defines a constant $p,1$-invariant graph.  For
$q \neq 1$ the situation is slightly more complicated, and there is no
direct analogue for $1,q$-invariant graphs in the unperturbed case.
The appropriate picture in this case is rather an invariant line for a
flow on the torus, as the following example shows. The fact that
invariant graphs may correspond to these two different types of
dynamics, or both at once, will come up again later, namely in the
proof of Thm.\ \ref{thm:denjoy}~.

Simple examples for $p,q$-invariant graphs with arbitrary $p$ and $q$
are given by torus translations: Consider $T : \ntorus \ra \ntorus, \ 
\thx \mapsto (\thom,x+\frac{k}{q}\omega + \frac{l}{pq})$ with $k$
relatively prime to $q$ and $l$ relatively prime to $p$.  Then
\[
\varphi^i_j(\theta) := \frac{k}{q}\theta + \frac{i-1+(j-1)p}{pq}
\hspace{4eM} (1 \leq i \leq p, \ 1 \leq j \leq q)
\] 
defines a continuous $p,q$-invariant graph. The meaning of the numbers
$k$ and $l$ will be explained in Section \ref{Rotnum}. Note that the
motion on $\Phi^i$ induced by the action of $T^p$ is equivalent to an
irrational rotation on \kreis\ (by
$\frac{p}{q}\omega+\frac{\tilde{l}}{q}$ where $l = \tilde{l}k \bmod
q$), such that it is not possible to decompose these graphs into
smaller invariant components.

\ \\
\textit{Furstenberg's results.}
As mentioned before, even in the absence of invariant graphs conjugacy
to an irrational torus translation cannot be expected, at least not if
we want the conjugacy to respect the skew product structure:
\begin{example} 
  A \textit{skew translation of the torus} is a map $T :\ntorus \ra
  \ntorus ,\ \thx \mapsto (\thom,x+g(\theta)\bmod 1)$ with some
  measurable function $g:\kreis \ra \R$. When $g$ is continuous, a
  fibrewise rotation number $\rho_T$ can be assigned to $T$ (see Def.\ 
  \ref{def:rotnum}), and $\rho_T$ will be equal to $\tau :=
  \int_{\kreis}g(\theta) \ d\theta$. As conjugacy preserves the
  rotation numbers, $T$ can only be conjugated to the irrational torus 
  translation $R:\thx \mapsto (\thom,x+\tau)$. 
  
  Now suppose $\tau \notin Q:= \{ \frac{k}{q}\omega+\frac{l}{p} \bmod
  1 \mid k,q,l,p \in \Z, \ p,q\neq 0\}$. Then $R$ is uniquely ergodic,
  and if $T$ and $R$ are conjugated (via a conjugacy $h$) so will be
  $T$. Thus $h$ must preserve the Lebesgue measure, which is invariant
  under both transformations. If we require the conjugacy to respect
  the skew product structure of the system, i.e.\ $\pi_1\circ h\thx =
  \theta \ \forall \thx \in \ntorus$, then all fibre maps $h_\theta$
  must also preserve Lebesgue measure and therefore be rotations.
  Hence $h$ will be a skew translation as well, i.e.\ $h\thx =
  (\theta,x+\varphi(\theta)\bmod 1)$ for some continuous function
  $\varphi:\kreis \ra \R$.  As $R=\hin \circ T\circ h$ this function
  $\varphi$ is a continuous solution of the cohomology equation
  \begin{equation} \label{eq:cohomology}
       g(\theta)-\tau = \varphi(\thom)-\varphi(\theta) \ .
  \end{equation}
  Conversely, if a function $g$ has measurable, but no continuous
  solutions of (\ref{eq:cohomology}), then on the one hand the
  resulting system $T$ cannot be conjugated to an irrational torus
  translation. On the other hand $T$ will still not have any invariant
  graphs, because a (non-continuous) solution $\varphi$ of
  (\ref{eq:cohomology}) still allows to define an isomorphism $h$
  between $T$ and $R$ as above. Thus $T$ will still be uniquely
  ergodic with respect to the Lebesgue measure on \ntorus.  An example
  of such a function $g$ is explicitly constructed in
  \cite{katok/hasselblatt:1997} (Section 12.6(b)).
\end{example}
This means that a possible analogue to the Poincar\'{e} classification
for circle homeomorphisms must be somewhat weaker in nature. It turned
out that the right perspective is to look at the invariant ergodic
measures (see \cite{furstenberg:1961}). First note that to any
invariant graph an invariant ergodic measure can be assigned:

\begin{definition}[Associated measure]
  \label{def:graphmeasure}
  Let $\varphi$ be a $p,q$-invariant graph. Then 
  \begin{equation}
    \label{graphmeasure}
    \muphi(A):= \frac{1}{pq} \sum_{i=1}^{p}\sum_{j=1}^q  m(\{\theta
    \in \kreis : 
    (\theta,\varphi^i_j(\theta)) \in A \} ) \ \ \forall A \in {\cal
      B}(\ntorus)  
  \end{equation}
  defines an invariant ergodic measure. If $\mu = \muphi$ for some
  $p,q$-invariant graph $\varphi$, then $\mu$ and $\varphi$ are called
  associated to each other.  
\end{definition}
In order to study general invariant ergodic measures, it is useful to
look at the so-called \textit{fibre measures}.  Any invariant measure
$\mu$ can be disintegrated in the way $\mu = m \times K$ where $K$ is
a stochastic kernel from \kreis\ to \kreis, such that $\mu(A) =
\int_{\kreis} K(\theta,\Ath) \ d\theta \ \forall A \in {\cal
  B}(\ntorus)$.  The measures $\muth := K(\theta,\ldot)$ on \kreis\ 
are the fibre measures of $\mu$. The fibre measures are mapped to each
other by the action of $T$, i.e.\ 
\begin{equation}
  \label{fibremeasure} 
  \muthom = \muth \circ T_{\thom}^{-1} \ \ \ \textrm{for }
  m\textrm{-a.e. } \theta \in \kreis \ .  
\end{equation}
All this is throughouly discussed for general random
dynamical systems in Chapter 1.4 in \cite{arnold:1998}. Note that an
invariant ergodic measure $\mu$ is associated to some invariant graph
if and only if its fibre measures are $m$-a.s.\ point
measures. Otherwise the fibre measures are continuous.

The following result of Furstenberg (\cite{furstenberg:1961},
Thm.\ 4.2) now provides a partial converse of Def.\ 
\ref{def:graphmeasure} and a substitute for the Poincar\'{e}
classification in the unforced case. Furstenberg used quite different
terminology, but his result can be stated in the following way:

\begin{thm}
  \label{thm:furstenberg}
  Let $T \in \Thom$. Then one of the following is true:
  \begin{list}{(\roman{enumi})}{\usecounter{enumi}}{}
  \item There exists a $p,q$-invariant graph $\varphi$. Then every
    invariant ergodic measure is associated to some $p,q$-invariant
    graph.
  \item There exists no invariant graph. Then $T$ is uniquely ergodic
    and the fibre measures of the unique invariant measure $\mu$ are
    continuous. 
\end{list}
\end{thm}
Indeed the original result is even slightly more general: 
Furstenberg only assumes that the base of the skew product is uniquely 
ergodic, not necessarily an irrational rotation. The main ingredient in 
the proof is the construction of an isomorphism $h$ between the
original system and a skew translation of the torus. This isomorphism
has some additional properties, which we describe by the following
\begin{definition}[Fibrewise conjugacy]
  \label{def:fibreconj}
Let $T,S \in \Thom$. $T$ is said to be fibrewise semi-conjugated to
$S$, if there exists a measurable map $h : \ntorus \ra \ntorus$ with the
following properties:
\romanlist
\item $\pie \circ h\thx = \theta \ \forall \thx \in \ntorus$
\item For every $\theta \in \kreis$ the map $h_\theta : \kreis \ra
  \kreis, \ x \mapsto \piz \circ h(\theta,x)$ is  continuous and
  order-preserving. 
\item $h \circ T = S \circ h$
\listend
If the fibre maps $h_\theta$ are all homeomorphisms, then $T$ and $S$
are called fibrewise conjugated. $h$ is called a fibrewise
semi-conjugacy or fibrewise conjugacy, respectively.
\end{definition}
Now Furstenbergs proof of Thm.\ \ref{thm:furstenberg} already contains
\begin{lem}
  \label{lem:isomorphism}
  Let $T \in \Thom$ and suppose there exists no invariant graph for
  $T$. Then $T$ is fibrewise semi-conjugated to an uniquely ergodic
  skew translation of the torus, and the semi-conjugacy $h$ maps the
  unique $T$-invariant measure $\mu$ to the Lebesque-measure on
  \ntorus.
\end{lem}
There is an interesting consequence of Thm.\ \ref{thm:furstenberg}
concerning Lyapunov exponents, which are defined pointwise by
$\lambda\thx := \nLim \log D\Tth^n(x)$. If the system is uniquely
ergodic all these limits exist and the convergence is uniform on
\ntorus.  As no iterate of $T$ can be uniformly expanding or
contracting, all Lyapunov exponents must be zero in the uniquely
ergodic case. Conversely, the existence of points with non-existent or
non-zero Lyapunov exponents implies that the system is not uniquely
ergodic and therefore has invariant graphs by Thm.\ 
\ref{thm:furstenberg} .

%
%
%                 Rotation numbers 
%
%

\section{Rotation numbers and regular invariant graphs} \label{Rotnum}

Throughout this and the next section, we will repeatedly work with
lifts of different objects on \ntorus\ to either $\kreis \times \R$ or
$\R^2$ (both of which are covering spaces of \ntorus). Therefore, we
have to make some conventions regarding notation and terminology:
\begin{itemize}
\item Projections: $\pi$ will denote the natural projection either from
  $\R^2 \ra \ntorus$, $\kreis \times \R \ra \ntorus$ or $\R \ra
  \kreis$. If we want to project to one coordinate and in addition
  project this coordinate onto the circle (if it is not in \kreis\ 
  anyway) we will denote this by $\pi_1$ and $\pi_2$. In those cases
  where we want to project to an \R-valued coordinate we will use
  $\widehat{\pi}_1$ or $\widehat{\pi}_2$, respectively. \R-valued
  variables will be denoted by $\hat{\theta},\hat{x},\hat{y},\ldots$ .
  The only exception is the rotation number $\omega$ on the base,
  which we will always identify with its unique lift in $[0,1)$.
\item A lift of a map $T \in \Thom$ to either $\kreis \times \R$ or
  $\R^2$ will usually be denoted by $\That$ (i.e.\ $T \circ \pi =
  \pi \circ \That$), for a second lift of the
  same map we will use $\Tbar$. It is  easy  to see that 
  \begin{equation}
    \label{eq:lifts}
    \That \ = \ \Tbar + (l,m) 
  \end{equation}
  for some $l,m \in \Z$ whenever $\That,\Tbar$ are two lifts of the
  same map $T$, and that a lift is uniquely determined by choosing its
  value $\That(\thhat_0,\xhat_0) = (\hat{\theta}_1,\hat{x}_1)$ at one
  point $(\thhat_0,\xhat_0)$. 
\end{itemize}
\textit{Fibrewise rotation numbers.} The fibrewise rotation number of
a quasiperiodically forced circle homeomophism can be defined exactly
in the same way as the rotation number in the unforced case:

\begin{thmdef}[Fibrewise rotation number]
  \label{def:rotnum}
  Let $T \in \Thom$  and $\That : \kreis \times  \R \ra  \kreis \times
  \R$ be a lift of $T$. Then the limit
  \begin{equation}
  \label{eq:rotnum}
        \rho_T := \nLim \ntel(\Thatth^n(\xhat) - \xhat) \ \bmod 1 \ 
  \end{equation}
  is independent of $\theta$, \xhat\ and the choice of the lift $\That :
  \kreis \times \R \ra \kreis \times \R$. It is called the fibrewise
  rotation number of $T$. Further more,
  \begin{equation}
  \label{eq:rotnumII}
  \rho_{\hat{T}} = \nLim \ntel \int_{\kreis} \widehat{T}_\theta^n(0) \
  d\theta 
  \end{equation}
  and the convergence in (\ref{eq:rotnum}) is uniform on $\kreis
  \times \R$.
\end{thmdef}
This result is due to Herman (\cite{herman:1983}), a very nice and
more elementary proof can be found in
\cite{stark/feudel/glendinning/pikovsky:2002}.  As we will see, the
existence of a continuous invariant graph implies that the fibrewise
rotation number \textit{depends rationally} on the rotation number on
the base, i.e.\
\[ 
    \rho_T  \ = \ \frac{k}{q}\omega  +  \frac{l}{n} \ \bmod 1
\]
holds for some $k,q,l,n \in \Z$. $\rho_T$ is called \textit{rationally
  independent} of $\omega$, if it is not rationally dependent.

\ \\
\textit{Dynamics of continuous invariant graphs.} In order to describe
the structure of continuous invariant graphs, we need the following
concept: Suppose a continuous function $\gamhat : \R \ra \R$ satisfies
\begin{equation}\label{eq:qcurve}
  \gamhat(\thhat + q) - \gamhat(\thhat) = k \ \textrm{ and } \ 
  \gamhat(\thhat + l) - \gamhat(\thhat) \notin \Z \ \ \forall \thhat \in
  \R,\ 1 \leq l < q
\end{equation}
for some $q \in \N$ and $k \in \Z$. We can then project it down to
\ntorus\ to obtain a $q$-valued graph
\[
  \gamma(\theta) = (\gamma_1(\theta) \ld \gamma_q(\theta)), \ \textrm{
    where } \ \gamma_i(\theta) = \pi(\gamhat(\thhat + i -1)) \hspace{3eM}
  (\thhat \in [0,1)),
\] 
which wraps around the circle $q$ times in the $\theta$-direction and
never intersects itself. In this situation we will call both $\gamma$
and the corresponding point set 
$\Gamma := \{ (\theta,\gamma_i(\theta)) \mid \theta \in \kreis, 1 \leq
i \leq q \}$ a \textit{$q$-curve}, a lift of $\gamma$ will always be a
function $\gamhat : \R \ra \R$ as in (\ref{eq:qcurve}). If $\gamma$ is
one-valued, it will be called a \textit{simple} curve.
  
Such a $q$-curve $\gamma$ has the property that it is a continuous
function $\kreis \ra \ntorus[q]_*$ and consists of one connected
component only (as a point set). Conversely, it follows from the usual
arguments for the existence of lifts that any graph $\gamma : \kreis
\ra \ntorus[q]_*$ with these properties is a $q$-curve, as it can be
lifted to a function $\gamhat : \R \ra \R$ which satisfies
(\ref{eq:qcurve}) and projects down to $\gamma$. Such a lift is
uniquely determined by choosing the value at one point, and by
$\gamhat_{x}$ we will denote the lift of $\gamma$ which satisfies
$\gamhat_x (0) = x$ (with $x \in \piin(\gamma(0))$).
  
Suppose $\varphi$ is a continuous $1,q$-invariant graph. Then $\Phi$
consists of one connected component only (otherwise the connected
components would be permuted, and $\Phi$ could therefore be decomposed
into periodic subgraphs). It follows from the comments made above that
$\varphi$ is a $q$-curve.  Similarly, any continuous $p,q$-invariant
graph $\varphi$ is a collection of $p$ disjoint $q$-curves $\varphi^1
\ld \varphi^p$.

\ \\
We now want to describe, how the structure and dynamics of an
invariant graph determine the rotation number. To that end, we need to 
define two characteristic numbers which can be assigned to an
invariant graph, the \textit{winding} and the \textit{jumping} number.

\begin{definition}[Winding number]
  \label{def:windnum}
Let $\gamma$ be a $q$-curve and $k \in \Z$, such that $\gamhat(\thhat
+ q) - \gamhat(\thhat) = k$ (see (\ref{eq:qcurve})). Then $k$ is
called the winding number of $\gamma$.
\end{definition}
In other words, the winding number is the number of times $\gamma$
winds around the circle in the $x$-direction before it closes. Some
simple observations about the winding number are collected in the
following:
\begin{lem}
  \label{lem:windnum}
Let $\gamma$ be a $q$-curve with winding number $k$. 
\alphlist
\item If $q >1$ then $k \neq 0$.

\item $q$ and $k$ are relatively prime.
\item If $\zeta$ is another  $q$-curve which does not intersect
  $\gamma$, then $\zeta$ has the same winding number as $\gamma$. 
\listend
\end{lem}
\proof \\
If $\gamhat$ is a lift of $\gamma$, then $\gamhat(\thhat + q) -
\gamhat(\thhat) = k \ \forall \thhat \in \R$. Thus, for any $j \in \N$
we have
\begin{eqnarray} 
    \lefteqn{\qtel \int_0^q \gamhat(\thhat + j) - \gamhat(\thhat) \
      d\thhat \  = } \nonumber  \\ & = & 
     \frac{j}{q} \int_0^q \gamhat(\thhat + 1) - \gamhat(\thhat) \ 
      d\thhat 
     \ = \ \frac{j}{q^2} \int_0^q \gamhat(\thhat + q) -
      \gamhat(\thhat) \ d\thhat \ = \ \frac{jk}{q} \label{eq:windnum}
\end{eqnarray}
\alphlist
\item Suppose $q > 1$ and $k = 0$. Taking $j = 1$ in
  (\ref{eq:windnum}) the Mean Value Theorem yields the existence of at
  least one $\thhat \in [0,q)$ with $\gamhat(\thhat + 1) =
  \gamhat(\thhat)$, a contradiction to the fact that $\gamma$ does not
  intersect itself (see (\ref{eq:qcurve})).
\item Suppose $q$ and $k$ are not relatively prime and take $j < q$,
  such that $\frac{jk}{q} \in \Z$. Then, as above, there is at least
  one $\thhat \in [0,q)$ which satisfies $\gamhat(\thhat + j) -
  \gamhat(\thhat) = \frac{jk}{q} \in \Z$, again contradicting
  (\ref{eq:qcurve}).
\item Suppose $\zeta$ has winding number $\tilde{k} < k$ and take
  lifts $\gamhat$ of $\gamma$ and $\widehat{\zeta}$ of $\zeta$ with
  $\widehat{\zeta}(0) \in [\gamhat(0), \gamhat(0)+1)$. But then
  $\widehat{\zeta}(q) = \widehat{\zeta}(0) + \tilde{k} < \gamhat(0) +
  k = \gamhat(q)$, and the two curves have to intersect somewhere in
  between. The case $\tilde{k} > k$ is treated analogously.
\listend 

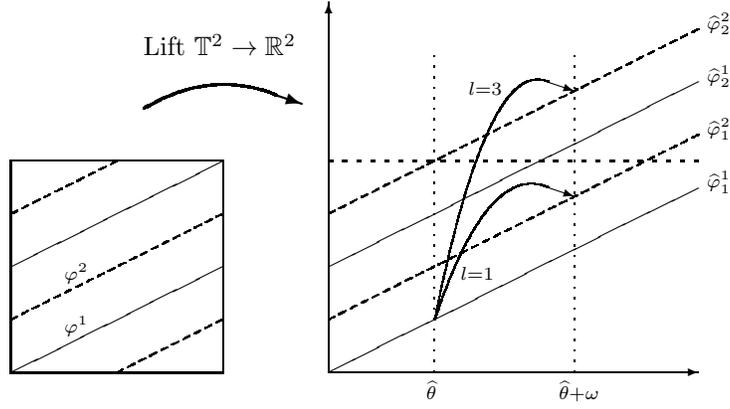
\begin{figure}[h!]
\begin{picture}(400,160)(-70,-10)
%
%                       Bildüberschrift
%
\put(10,180){}

\put(0,0){\framebox(80,80)}
\drawline(0,0)(80,40)
\drawline(0,40)(80,80)

\dashline[+50]{3}(0,20)(80,60)
\dashline[+50]{3}(0,60)(40,80)
\dashline[+50]{3}(40,0)(80,20)

\put(20,15){$\scriptstyle \varphi^1$}
\put(20,35){$\scriptstyle \varphi^2$}

\put(120,0){\vector(0,1){140}}
\put(120,0){\vector(1,0){140}}

\dashline{2}(120,80)(260,80)

\drawline(120,0)(260,70)
\drawline(120,40)(260,110)

\dashline[+50]{3}(120,20)(260,90)
\dashline[+50]{3}(120,60)(260,130)

\put(263,70){$\scriptstyle \phihat^1_1$}
\put(263,90){$\scriptstyle \phihat^2_1$}
\put(263,110){$\scriptstyle \phihat^1_2$}
\put(263,130){$\scriptstyle \phihat^2_2$}

\qbezier(160,20)(180,80)(203,70)
\qbezier(160,20)(180,120)(203,110)

\put(170,35){$\scriptstyle l=1$}
\put(173,105){$\scriptstyle l=3$}

\dashline[+50]{1}(160,0)(160,120)
\put(157,-10){$\scriptstyle \thhat$}

\dashline[+50]{1}(213,0)(213,120)
\put(207,-10){$\scriptstyle \thhat+\omega$}

\put(203,70){\vector(3,-1){10}}
\put(203,110){\vector(3,-1){10}}
\thicklines
\qbezier(50,100)(75,115)(100,105)
\put(100,105){\vector(3,-1){10}}
\put(50,120){Lift $\ntorus \ra \R^2$}
\end{picture}
\caption{\small Depicted on the left is the graph $\varphi =
  (\varphi^1,\varphi^2)$ with $\varphi^i_j(\theta) = \halb\theta +
  \frac{i-1+2(j-1)}{4} \ (i,j=1,2)$. This is a $2,2$-invariant graph for both $\thx
  \mapsto (\thom,x+\halb\omega+\viertel)$ and $\thx \mapsto
  (\thom,x+\halb\omega+\dreiviertel)$, but the jumping number is
  different in the two cases ($l=1$ and $l=3$
  respectively). The different lifts mentioned before Def.\
  \ref{def:jumpnum} are shown on the right.}
\label{fig:liftedgraphs}
\end{figure}

\qed

\ \\
In order to define the jumping number, suppose $T \in \Thom$ has a
continuous $p,q$-invariant graph $\varphi = (\varphi_1^1 \ld
\varphi^p_q)$, consisting of $p$ $p$-periodic $q$-curves $\varphi^1 \ld
\varphi^p$. Let
$\xhat_1^1 \ld \xhat^p_q \in[0,1)$
be lifts of $\varphi_1^1(0) \ld \varphi^p_q(0)$. W.l.o.g.\ we can
assume that these points are ordered in the way
\[
\xhat^1_1 < \xhat^2_1 < \ldots < \xhat^p_1 < \xhat^1_2 < \ldots <
\xhat^p_2 < \ldots < \xhat^p_q \ .
\] 
(Note that any of the intervals $[\varphi^1_j(0),\varphi^1_{j+1}(0)]$
must contain the same number of points from the other graphs
$\varphi^i$.) To each of the points $\xhat^i_j$ we can assign a lift
$\phihat^i_j := \phihat^i_{\xhat^i_j}$ of the $q$-curve $\varphi^i$,
and there is exactly one lift $\That$ of $T$ which maps $\phihat^1_1$
to one of the other lifts $\phihat^i_j$ (see Fig.\
\ref{fig:liftedgraphs} for an example).

\begin{definition}[Jumping number]
  \label{def:jumpnum}
  Let $0 \leq m \leq p -1$ and $0 \leq n \leq q -1$ be such, that
  there is a lift $\That$ of $T$ for which
  \[
  \That(\thhat,\phihat^1_1(\thhat)) =
  (\thhat+\omega,\phihat^{1+m}_{1+n}(\thhat+\omega)) \ .
  \]
  Then $l := m + np$ is called the jumping number of $\varphi$ (with
  respect to $T$).
\end{definition}

\begin{bem}
  \label{bem:jumpnum}
  If $l = m + np$ is the jumping number of a $p,q$-invariant graph
  $\varphi$, then $m$ and $p$ are relatively prime, which is
  equivalent to $l$ and $p$ being relatively prime. Otherwise $p \mid
  m p'$ holds for some $p' \in \{1 \ld p-1\}$. But then the union of
  $\Phi^1, \Phi^{1+m} \ld \Phi^{1+(p'-1)m}$ is invariant,
  contradicting the minimality of $\Phi$.
\end{bem}

\begin{prop}
  \label{prop:rotnum}
Suppose $T \in\Thom$ and $\varphi$ is a continuous $p,q$-invariant
graph with winding number $k$ and jumping number $l$. Then 
\begin{equation}
  \label{eq:ctsrotnum}
  \rho_T \ = \ \frac{k}{q}\omega + \frac{l}{pq} \ \bmod 1 \ .
\end{equation}
Conversely, the numbers $p,q,k$ and $l$ associated to a
continuous invariant graph in the above way are uniquely determined by
the rotation numbers $\omega$ and $\rho_T$.
\end{prop}
\proof \\
Obviously, instead of working with a lift to $\kreis \times \R$ we can
also work with a lift to $\R^2$. By Def./ Thm.\ \ref{def:rotnum}
we can use any such lift, any initial point $(\thhat,\xhat)$ and any
sequence of iterates to determine the fibrewise rotation number. Let
therefore $\phihat^1_1 \ld \phihat^p_q$ be as in Def.\ 
\ref{def:jumpnum} and choose the lift
$\That : \R^2 \ra \R^2$, such that
 $\That(0,\phihat^1_1(0)) = (\omega,\phihat^{1+m}_{1+n}(\omega))$.
 After $pq$ iterates, we will return to a copy of $\phihat^1_1$ which
 is vertically translated exactly by $l$, i.e.\ 
\[  
    \Thatth^{pq}(\phihat^1_1(0)) = \phihat^1_1(pq\omega) + l \ .    
\] 
Thus we get
\[
\rho_{\That} = \nLim \frac{1}{npq} \left(\phihat^1_1(npq\omega) -
  \phihat^1_1(0) \right) + \frac{l}{pq}\ .
\]
By using the continuity of $\phihat^1_1$ and (\ref{eq:qcurve}) we can
replace the irrational translation with integer translation, i.e.\ 
\[
\nLim \frac{1}{npq}\left(\phihat^1_1(npq\omega) -
  \phihat^1_1(0)\right) = \lim_{\xhat\ra \infty} \frac{1}{\xhat}
\left(\phihat^1_1(\xhat \omega) - \phihat_1^1(0) \right) = \jLim
\frac{\omega}{jq} (\phihat^1_1(jq) - \phihat^1_1(0)) =
\frac{k}{q}\omega \ .
\]
Putting all this together yields $\rho_{\That} = \frac{k}{q} \omega +
\frac{l}{pq}$.
\ \\
It remains to show, that $p,q,k$ and $l$ are uniqely
determined by $\omega$ and $\rho_T$. To that end, suppose
\[
    \frac{k}{q} \omega + \frac{l}{pq}   =   \frac{k'}{q'} \omega +
    \frac{l'}{p'q'}  \ ,
\]
and the numbers on both sides belong to continuous invariant graphs.
Then by Lemma \ref{lem:windnum} $k$ and $q$ as well as $k'$ and $q'$
must be relatively prime and therefore $k = k'$ and $q = q'$. The same
argument then applies to $l,p$ and $l',p'$ (see Rem.\ 
\ref{bem:jumpnum}).

\qed

\ \\
\textit{Regular invariant graphs.}  So far we have seen that the
existence of continuous invariant graphs allows to draw conclusions
about the fibrewise rotation number. This stays true if the assumption
of continuity is weakened to some extent, namely for the class of
\textit{regular} invariant graphs which we will introduce in this
section. This concept will then turn out to be vitally important for
both the formulation and the proof of Thm.\ \ref{thm:denjoy}~.  The
regularity we require is, that the invariant graphs are boundary lines
of certain compact invariant sets. This is hardly surprising,
considering that such boundary lines, which are necessarily
semi-continuous functions, also play a prominent role in the study of
quasiperiodically forced interval maps (see e.g.\ 
\cite{keller:1996,stark:2003}). However, instead of working with the
compact invariant sets it will be more convenient for our purposes to
use their complements.

\begin{definition}[$p,q$-invariant strips and $p,q$-invariant open tubes]
  \label{def:tubes}
  Let $T \in \Thom$.
\alphlist
\item A compact invariant set $A$ is called a $1,q$-invariant strip if
  it consists of $q$ disjoint closed and nonempty intervals on every
  fibre and only has one connected component. When $A^1$ is a
  $1,q$-invariant strip for $T^p$ and the sets $A^i := T^{i-1}(A^1) \ 
  (i=2\ld p)$ are pairwise disjoint, then $A:= \bigcup_{i=1}^{p} A^i$
  is called a $p,q$-invariant strip.
\item Similarly, an open invariant set $U$ is called a $1,q$-invariant
  open tube, if it consists of $q$ disjoint open and nonempty
  intervals on every fibre and only has one connected component. When
  When $U^1$ is a $1,q$-invariant open tube for $T^p$ and the sets $U^i :=
  T^{i-1}(U^1) \ (i=2\ld p)$ are pairwise disjoint, then $U:=
  \bigcup_{i=1}^{p} U^i$ is called a $p,q$-invariant strip.
\listend
\end{definition}
Of course a continuous invariant graph is a special case of an
invariant strip. By some elementary topological arguments one can also 
see that the above definitions are indeed complementary:
\begin{bem}
  \label{bem:tubes}
  \alphlist
\item The complement of a $p,q$-invariant strip is always a
  $p,q$-invariant open tube and vice versa. 
\item If $U = \bigcup_{i=1}^p U^i$ is a $p,q$-invariant open tube,
  then any of the connected components $U^i$ contains a $q$-curve
  $\gamma^i$ such that each of the $q$ intervals of $U^i_\theta$
  contains exactly one of the points $\gamma^i_1(\theta) \ld
  \gamma^i_q(\theta)$. As the $q$-curves $\gamma^1 \ld \gamma^p$
  cannot intersect, they all have the same winding number
\item A winding and a jumping number can be assigned to
  invariant open tubes in almost exactly the same way as to
  continuous invariant graphs: The winding number of a $q,p$-invariant 
  open tube $U$ is that of the $q$-curves it contains, and the jumping
  number can be defined  compleatly analogous to Def.\
  \ref{def:jumpnum}~ by looking at the different lifts of the
  $U^i$ instead of those of the $\varphi^i$.
\listend
\end{bem}
\begin{lem}
  \label{lem:tubes}
Let $T \in \Thom$ and suppose there exists a $p,q$-invariant open tube 
with winding number $k$ and jumping number $l$. Then 
\[
     \rho_{T} \ = \ \frac{k}{q} \omega + \frac{l}{pq} \ \bmod 1 \ .
\]
Conversely, the numbers $p,q,k$ and $l$ are uniquely determined by
$\omega$ and $\rho_T$.
\end{lem}
\proof \\
The proof is almost the same as for Prop.\ \ref{prop:rotnum}.
Instead of the lift $\phihat^1_1$ of $\varphi^1$ we use a lift
$\gamhat^1_1$ of the $q$-curve $\gamma^1$. Then $U_i$ can be
lifted around $\gamhat_1^1$ in a unique way, and although $\gamma$ is
not invariant the properties of $U_i$ still guarantee
\[
    \left| \Thatth^{npq}(\gamhat^1_1(0)) - \gamhat^1_1(npq\omega) -
      nl \right| \ < \ 1 \ .
\]
From that the result follows easily by repeating the proof for
continuous invariant graphs, replacing $\phihat^1_1$ with $\gamhat^1_1$. 

\qed

\ \\
Invariant strips and tubes are of course bounded by invariant graphs,
which are semi-continuous in a sense (although semi-continuity is hard 
to define for \kreis-valued functions). These invariant graphs could
be defined by just taking the right (or left) endpoints of the
intervals on the fibres as values, but with the $q$-curves contained
in an invariant open tube this can be made a little bit more precise:

\begin{definition}[Regular invariant graphs]
 \label{def:regularinvgraphs}
  $\varphi : \kreis \ra \ntorus[pq]_*$ is called a regular
  $p,q$-invariant graph if it is the boundary line of a
  $p,q$-invariant open tube, i.e.\ 
\[
    \varphi^i_j(\theta) = \gamma^i_j(\theta) + \inf \{x \in [0,1) \mid 
    \gamma^i_j(\theta) + x \in U_\theta^c \}
\]
or
\[
\varphi^i_j(\theta) = \gamma^i_j(\theta) - \sup \{ x \in [0,1) \mid
\gamma^i_j(\theta) - x \in U_\theta^c \} \ ,
\]
where $U = \bigcup_{i=1}^{pq} U^i$ and the $\gamma^i$ are $q$-curves
contained in $U^i$. 
\end{definition}

%
%
%                              Denjoy
%
%

\section{The Denjoy argument} \label{Denjoy}

The aim of this section is to carry over the distortion argument used
in the proof of the Theorem of Denjoy to the quasiperiodically forced
case.  Denjoy's Theorem states, that a circle diffeomorphism with
irrational rotation number and bounded variation of the derivative is
not only semi-conjugated but conjugated to the corresponding rotation.
As a simple consequence of this the system will be minimal.
However, the crucial step in the proof is to exclude the possibility
of wandering open sets, and the rest then follows quite easily from
the fact that a semi-conjugacy to the irrational rotation is already
established. This is different in the quasiperiodically forced case,
where the absence of invariant graphs does not imply semi-conjugacy to 
an irrational rotation. Thus the Denjoy argument generalized to the
forced case will only yield that there are no wandering open sets. A
few simple arguments can then be used to establish transitivity (see
Thm.\ \ref{thm:denjoy}), but whether such systems are minimal has to be
left open here.

There are two other problems which make the quasiperiodically forced
case more complicated: On the one hand, in the one-dimensional case
the semi-conjugacy result provides all the combinatorics which are
necessary. Without such a result in the forced case this has to be
done ``bare hands'', and thus a significant part of the proof of Thm.\ 
\ref{thm:denjoy} will consist of dealing with the possible
combinatorics of wandering open sets.

Another fact was already mentioned in Section \ref{Furstenberg}: An
invariant graph may correspond to two different dynamical situations.
Thus the distortion argument will have to be applied twice, first in
the proof of Lemma \ref{lem:closestret} and then in the proof of Thm.\ 
\ref{thm:denjoy}, each time dealing with one of the two possibilities.
To get a better picture of this, one should recall Denjoy's examples
of diffeomorphisms which are not conjugate to the corresponding
irrational rotation. The construction of these starts with an orbit of
the irrational rotation, which is then ``blown up'' to yield wandering
intervals. Now, in the quasiperiodically forced case there are two
objects which might be used for such a construction. One is just the
orbit of a constant line, similar to the one-dimensional case, but
instead one could also use an infinite invariant line for a flow on
the torus.

\ \\
If the fibre maps $\Tth$ of $T \in \Thom$ are differentiable for
$m$-a.e.\ $\theta$, let 
\[
    V(T) \ := \ \int_{\kreis} V_\theta \ d\theta \ ,
\]
where $V_\theta$ denotes the variation of $\log D\Tth$, i.e.\ 
$V_\theta = \sup \{ \insum |\log D\Tth(x_i) - \log D\Tth(x_{i-1})| \ 
\mid 0 \leq \xeins \leq \ldots \leq \xn \leq 0\}$. By ${\cal T}_{BV}$ we
will denote the set of all such $T$ with $V(T) < \infty$. 

Bounded variation usually implies bounded distortion of the derivative
(Compare Lemma 12.1.3 in \cite{katok/hasselblatt:1997} or Corally 2 to
Lemma 2.1 in \cite{demelo/vanstrien:1993}). In the forced case,
however, any such statement must always be an integrated version of
the one-dimensional. To obtain fibrewise bounds will hardly be
possible, as an equal amount of distortion might be picked up from a
different fibre in each iteration step. For the formulation of the
lemma below we will also consider graphs $\varphi,\psi : I \ra \kreis$
which are defined only on a subinterval $I \subseteq \kreis$ and use
the notation
\[
    [\varphi,\psi] \ := \ \{ \thx \in \ntorus \mid \theta \in I, \ x
    \in [\phith,\psith] \} \ .
\]
Then we have
\begin{lem}
  \label{lem:distortion}
Let $T \in \Tbv$, $I \subseteq \kreis$ an interval and $\varphi,\psi
: I \ra \kreis$ be such that $[\varphi,\psi],T([\varphi,\psi]) \ld
T^{n-1}([\varphi,\psi])$ are pairwise disjoint. Then $\forall s > 0$
\[
     \int_I \left( \frac{D\Tnth(\psith)}{D\Tnth(\phith)}\right)^s \
     d\theta \ \geq \ |I| e^{-\frac{sV(T)}{|I|}} \ .
\]
\end{lem}
\proof \\
\begin{eqnarray*}
  \lefteqn{ \log \int_I \left(
      \frac{D\Tnth(\psith)}{D\Tnth(\phith)}\right)^s \ 
     \cdot \frac{1}{|I|} \ d\theta  \ \stackrel{\scriptstyle
       \textrm{Jensen}}{\geq} 
     \int_I \log \left(\frac{D\Tnth(\psith)}{D\Tnth(\phith)}\right)^s \
     \cdot \frac{1}{|I|} \ d\theta  \ = \ } \\
   & = &
     \frac{s}{|I|} \int_I \inergsum \left( \log
       DT_{\theta+i\omega}(\Tth^i(\psith) - \log DT_{\theta +
         i\omega}(\Tth^i(\phith) \right) \ d\theta \ = \ \\
   & = & 
     \frac{s}{|I|} \int_{\kreis} \underbrace{\twoscriptsum{0 \leq i <
         n}{\theta-i\omega \in I} \hspace{-1eM} \left( \log
       DT_{\theta}(T_{\theta-i\omega}^i(\psi(\theta-i\omega)) - \log
       DT_{\theta}(T_{\theta - i\omega}^i(\varphi(\theta-i\omega))
     \right)}_{\geq -V_\theta} \ d\theta \ 
   \geq \ - \frac{s}{|I|} V(T) \ \\ 
\end{eqnarray*}

\begin{figure}[h!] 
\noindent
\begin{minipage}[t]{\linewidth}
 \hspace{10eM} \epsfig{file=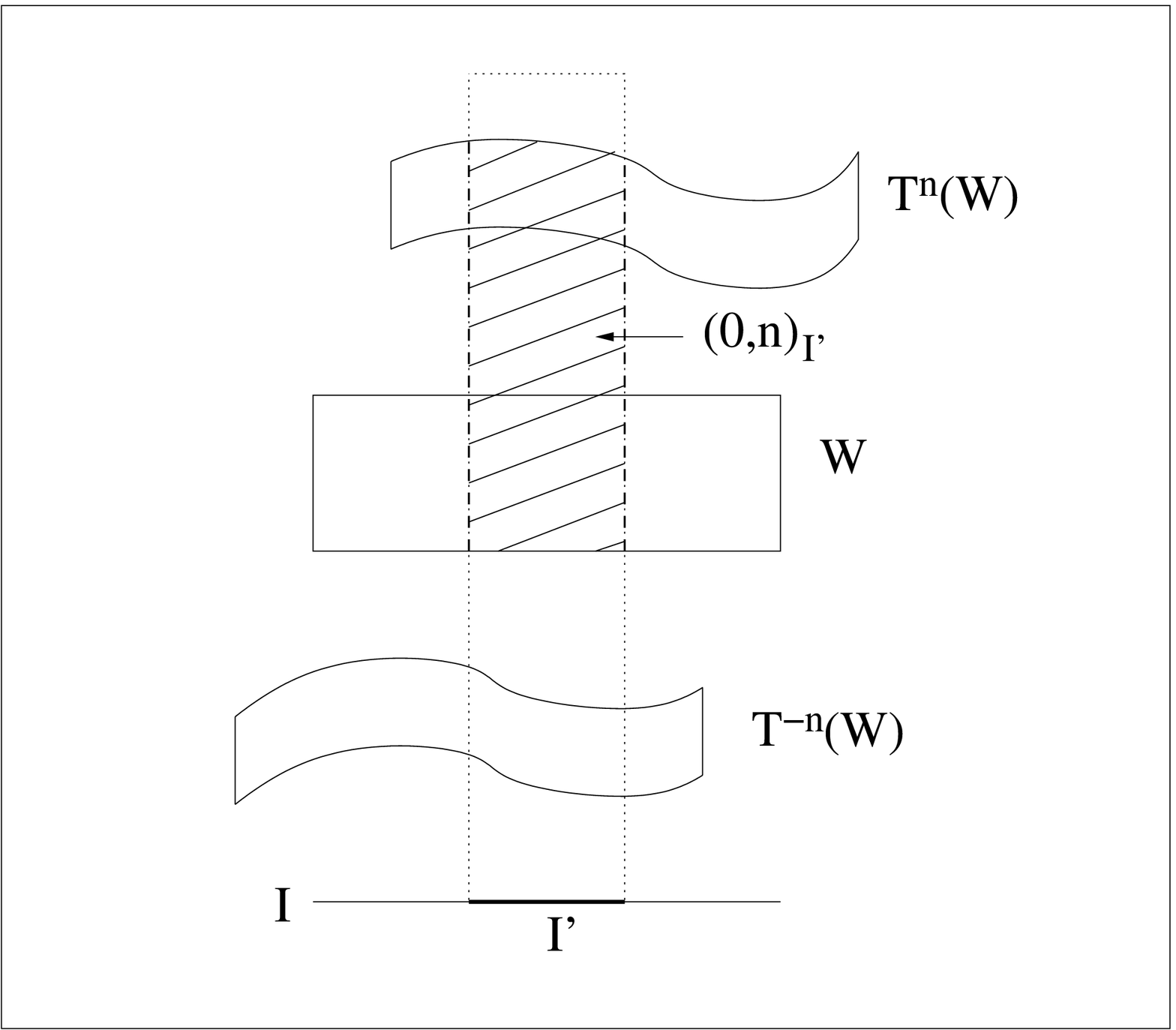, height=0.3\textheight,
  width=0.5\textwidth} \caption{Here $n$ is comparable over
  $I'=I_{\drittel}$ and $-n \lhd 0 \lhd n$. } \label{combinatorics}
\end{minipage}
\end{figure}

\qed

\ \\
In order to apply this lemma, we will now look at the
combinatorics of wandering boxes. To that end, some
notation must be introduced (see also Fig.\ \ref{combinatorics}): 
\begin{itemize}
\item In all of the following we will assume that $W = I \times K$ is
  a wandering set, i.e.\ $T^n(W) \cap W = \emptyset \ \forall n \in
  \N$, where $I,K \subset \kreis$ are open intervals and $|I| <
  \halb$.  For the sake of simplicity we will assume that $I$ is
  centered around $0$. $I_\alpha$ will denote the symmetric
  middle part of $I$ with length $\alpha|I|$, i.e.\ $I_\alpha :=
  (-\frac{\alpha|I|}{2},\frac{\alpha|I|}{2})$.
\item $\neins \ld \nk \in \N$ are called \textit{comparable} over an
  interval $J$, if $J \subseteq \ikcap I + n_i\omega$ \ . 
\item If $k\geq 3$, then 
   \[
   \neins \lhd \ldots \lhd \nk \lhd \neins \ \textrm{ over } J
   \]
   means that $\neins \ld \nk$ are comparable over $J$ and 
  \[
  \forall \theta \in J \ \forall x_i \in (T^{n_i}W)_\theta : \xeins <
  \ldots < x_k< \xeins
   \]
   (if this last statement is satisfied for some
   $\theta \in J$ and $x_i \in (T^nW)_\theta$, then it is true for
   all). We will usually omit the last $\lhd \neins$ and only write
   $\neins \lhd \ldots \lhd \nk$.
\item Let $n_1$ and $n_2$ be comparable over $J$. Then
   \[
        (n_1,n_2)_J \ := \ \{ \thx \in \ntorus \mid \theta \in J
        \textrm{ and } \exists x_1 \in (T^{n_1}W)_\theta,\ x_2 \in
        (T^{n_2}W)_\theta : x_1 \leq x \leq x_2 \} \ .
   \]
 \item The set of return times (with respect to $I_\alpha$) is defined
   as 
   \[
   N(\alpha) := \{ n \in \Z \mid n \textrm{ comparable over } I_\alpha
   \} = \{ n \in \Z \mid |n\omega \bmod 1| \leq \frac{1-\alpha}{2}|I|
   \} \ .
   \]
 \item $n \in N(\alpha)$ is called a \textit{closest return time}
   (with respect to $I_\alpha$), if \ 
  \[
  -n \lhd 0 \lhd n \ \textrm{ over } I_\theta \ \textrm{ and } \ 
  \nexists k \in N(\alpha)\setminus\{0\} : |k| < |n|, \ 0 \lhd k \lhd
  n \ \textrm{ over } I_\alpha
  \]
  or
  \[
  n \lhd 0 \lhd -n \ \textrm{ over } I_\theta \ \textrm{ and } \ 
  \nexists k \in N(\alpha)\setminus\{0\}: |k| < |n|, \ n \lhd k \lhd 0
  \ \textrm{ over } I_\alpha \ .
  \]
\end{itemize}
Note that it is usually not possible to use the ``ordering'' $\lhd$
just in a formal way, as the interval with respect to which it is used
always has to be taken into account. Therefore, a little bit of care
has to be taken and in the following remark a few simple facts are
collected. They might all seem trivial, but as they are used
frequently in different combinations later on this should help to
avoid confusion.

\begin{bem}
  \label{bem:closestret}
\romanlist
\item If $J$ is symmetric in $I$ (i.e.\ $J=I_\alpha$ for some
  $\alpha$) and $n,k$ are comparable over $J$,
  then either $k$ is comparable over $J + n \omega$ (this is the case
  when $n\omega \bmod 1$ and $k\omega \bmod 1$ are in the same half of
  $I$), or $k$ is comparable over $J - n\omega$ and $n+k$ is comparable
  over $J$ (when $n\omega \bmod 1$ and $k\omega \bmod 1$ are in
  opposite halfs of $I$ they ``cancel each other out''). Thus $0,n,k$
  and $n+k$ will always be comparable over one of the two intervals
  $J$ or $J+n\omega$.
\item $\neins \lhd \ldots \lhd \nk \ \textrm{ over } J \ \equi \ n_1 + 
  i \lhd \ldots \lhd \nk + i \ \textrm{ over } J + iw \ \ \forall i
  \in \Z$.
\item $\neins \lhd n_2 \lhd n_3$ \ over $J$ and $\neins \lhd n_3 \lhd
  n_4$ \ over $J$ \ \follows\ \ $\neins \lhd n_2 \lhd n_3 \lhd n_4$ \ 
  over $J$.
\item Let $n_1 \ld \nk$ be comparable both over $J_1$ and $J_2$. Then
  \[
  n_1 \lhd \ldots \lhd \nk \ \textrm{ over } J_1 \ \equi \ n_1 \lhd
  \ldots \lhd n_k \ \textrm{ over } J_2 \ .
  \]
\item $-n \lhd 0 \lhd k \lhd n$ \ over $I_\alpha$ \ $\equi$ \ $-n \lhd
  -k \lhd 0 \lhd n$ \ over $I_\alpha$. In particular, if $n$ is a
  closest return time then so is $-n$.  \listend
  \ \\
  Apart from the last all these statements should be obvious, and (v)
  can be seen as follows: Due to the symmetry of $I_\alpha$, $k$ is
  either comparable over $I_\alpha + n\omega$ or over $I_\alpha -
  n\omega$ (compare (i)). Suppose $k$ is comparable over
  $I_\alpha+n\omega$. Then by (ii) $-n \lhd k-n \lhd 0 \ $ \ over
  $I_\alpha$ and this can be extended to $-n \lhd k-n \lhd 0 \ \lhd k$
  \ over $I_\alpha$ by (iii). Thus $k-n \lhd 0 \lhd k$ \ over
  $I_\alpha + k\omega$ by (iv) and applying (ii) and (iii) for a
  second time yields $-n \lhd -k \lhd 0 \lhd n$ \ over $I_\alpha$.
  The case where $k$ is comparable over $I_\alpha-n\omega$ is treated
  similarly: Then $-n \lhd k-n \lhd 0 \lhd k$ \ over
  $I_\alpha-n\omega$ by (ii) and (iii) and therefore $k-n \lhd 0 \lhd
  k$ \ over $I_\alpha + k\omega$ by (iv) as before.
\end{bem}
With these notions it is now possible to give a lower bound for the
size of the images of $W$ at closest return times:
\begin{lem}
  \label{lem:closestret}
  Let $T \in \Tbv$ and suppose $n \in N$ is a closest return time with
  respect to $I_\alpha$, w.l.o.g.\ $-n \lhd 0 \lhd n$ \ over
  $I_\alpha$, $\beta := \min\{ \alpha, \frac{1-\alpha}{2} \}$.  Then
  \romanlist
\item   
$
    T^k(0,n)_{I_\beta} \cap (0,n)_{I_\beta} = \emptyset \ \ \forall k :
    |k| < |n|
$
\item 
$
     \lambda(T^nW \cup T^{-n}W) \ \geq \ |K|\cdot|I_\beta|\cdot
     e^{-\frac{V(T)}{2|I_\beta|}}  \ .
$
\listend
In particular, there can only be finitely many closest returns. 
\end{lem}
\proof 
\romanlist
\item Suppose for a contradiction that $T^k(0,n)_{I_\beta} \cap
  (0,n)_{I_\beta} \neq \emptyset$. Of course, this can only be if
  $I_\beta \cap (I_\beta + k\omega) \neq \emptyset$, i.e.\ $|k\omega
  \bmod 1| < \beta$, and this implies $k \in N(\alpha)$.  Further,
  there are in principle three possibilities: $0 \lhd k \lhd n$ or $0
  \lhd n+k \lhd n$ or $k \lhd 0 \lhd n \lhd n+k$. Of course we have to
  specify with respect to which intervals these numbers can be
  compared, in particular in the last case which has to be split up
  once more (according to Rem.\ 
  \ref{bem:closestret}(i)). \\
  Case 1: $0 \lhd k \lhd n$ \ over $I_\alpha$. \\
  Case 2: $0 \lhd n+k \lhd n$ \ over $I_\alpha + n\omega$. But then by
  Rem.\ \ref{bem:closestret}(ii) $-n \lhd k \lhd 0$ \ over $I_\alpha$. \\
  Case 3(a): $k \lhd 0 \lhd n \lhd n+k$ over $I_\alpha$. Then $0 \lhd
  -k \lhd n-k \lhd n$ \ over $I_\alpha -k\omega$ and thus also $0 \lhd
  -k \lhd n$ \ over $I_\alpha$ (Rem.\ 
  \ref{bem:closestret}(ii) and (iv) respectively). \\
  Case 3(b): $k \lhd 0 \lhd n \lhd n+k$\ over $I_\alpha+k\omega$. Then
  $0 \lhd -k \lhd n-k \lhd n$ \ over $I_\alpha$ (Rem.\ 
  \ref{bem:closestret}(ii)).
  \ \\
  In all cases this contradicts the fact that $n$ (and thus $-n$ by
  Rem.\ \ref{bem:closestret}(v)) is a closest return time.
\item Note that $I_\beta + n\omega \subseteq I$. Thus
  \begin{eqnarray*}
    \lefteqn{\lambda(T^nW \cup T^{-n}W) \ = \ \int_K \int_I D\Tnthx +
      D\Tth^{-n}(x) \ d\theta dx \ \geq } \\
    & \geq & \int_K \int_{I_\beta} D\Tnthx +
      DT_{\theta+n\omega}^{-n}(x) \ d\theta dx  \ \geq \ \int_K
      \int_{I_\beta}\left( D\Tnthx \cdot 
      DT_{\theta+n\omega}^{-n}(x)\right)^{\halb} \ d\theta dx \ = \\
    & = &
    \int_K \int_{I_\beta} \left(
      \frac{D\Tnthx}{D\Tnth(T^{-n}_{\theta+n\omega}(x))} 
    \right)^{\halb} d\theta dx \ \stackrel{\scriptstyle
      \textrm{Lem.\ \ref{lem:distortion}}}{\geq} \ \int_K |I_\beta| 
      e^{-\frac{V(T)}{2|I_\beta|}} \ dx \ = \ |K| |I_\beta|
      e^{-\frac{V(T)}{2|I_\beta|}} \ .
  \end{eqnarray*}
\listend
As $W$ is wandering we have that $\sum_{n \in \Z} \lambda(T^nW)
\leq 1$, and therefore the number of closest returns must be finite. 

\qed 

\ \\
This lemma will turn out to be crucial in the proof of the next
theorem, which is the main result of this section.

\begin{thm}
  \label{thm:denjoy}
  If $T \in \Tbv$ is not topologically transitive, then there exists a
  regular $p,q$-invariant graph for $T$. In particular, $\rho_T$
  depends rationally on $\omega$.
\end{thm}
\proof \\
The proof consists of two steps. The first is to show that the
existence of a wandering open set implies the existence of a
regular invariant graph. The second will then be to show, that the
non-existence of wandering sets already implies the transitivity of
the system. For this, we will use the fact that no iterate of $T$ can
have wandering open sets in the absence of regular invariant graphs
either, as all these iterates are also in $\Tbv$ and therefore Step 1
applies to them as well.

\ \\
\underline{Step 1:} \ By Lemma \ref{lem:closestret} there can only be
finitely many closest return times with respect to $I_{\halb}$. Let
$p$ be the maximum of these, w.l.o.g.\ $-p \lhd 0 \lhd p$ \ over
$I_{\halb}$.  As a consequence of Rem.\ \ref{bem:closestret}(v) we
have $\nexists n \in \Z \setminus \{ 0 \}: -p \lhd n \lhd p$ \ over
$I_{\halb}$, otherwise the minimum of such $n > 0$ would be a closest
return time greater than $p$. Of course 
$T^p(0,p)_{I_{\halb}} \cap (0,p)_{I_{\halb}} \neq \emptyset$, but as
in the proof of Lemma \ref{lem:closestret}(i) we can conclude that
afterwards
\begin{equation}
  \label{eq:lastret}
  T^{np}(0,p)_{I_{\halb}} \cap (0,p)_{I_{\halb}} \ = \ \emptyset \ \ \
  \forall n \geq    2 \ .
\end{equation}
It follows that $\bigcup_{n\in\Z} T^{np}(0,p)_{I_{\halb}}$ is an
$T^p$-invariant open set and looks like a ``tube'' which winds around the
circle infinitely many times in the $\theta$-direction. To make this
more precise, we will construct a infinite invariant line $\gamma$
inside of this set which never intersects itself.

\ \\
Let $\omtil := p\omega \bmod 1$. For the sake of simplicity we assume
that $\omtil$ is close to $0$ from the right, e.g.\ $\omtil \in
[0,\viertel)$. The open set $(0,p)_{I_{\halb}} \cup
T^p(0,p)_{I_{\halb}}$ is connected, thus there exists a continuous
function $\gamma_0 : [0,\omtil] \ra \kreis$, such that
$(0,\gamma_0(0)) \in (0,p)_{I_{\halb}}$, $\gamma_0(\omtil) =
T_0^p(\gamma_0(0))$ and $(\theta,\gamma_0(\theta)) \in
(0,p)_{I_{\halb}} \cup T^p(0,p)_{I_{\halb}} \ \forall \theta \in
[0,\omtil]$.  Let $\gamhat_0 : [0,\omtil] \ra \R$ be a lift of
$\gamma_0$, $\That : \R^2 \ra \R^2$ a lift of $T^p$. Then by
\[
    \gamhat(\thhat) := \That_{\theta-n\omtil}^n(\gamhat_0(\thhat -
    n\omtil) ) \ \ \textrm{ if } \ \thhat \in [n\omtil,(n+1)\omtil)
\]
a continuous $\That$-invariant curve $\gamhat : \R \ra \R$ can be
defined. $\gamhat$ projects down to an infinite $T^p$-invariant curve
$\gamma : \R \ra \kreis$, and the point set $\Gamma := \{
\pi(\thhat,\gamma(\thhat)) \mid \thhat \in \R \} \subseteq \ntorus$
lies inside of $\cup_{n\in\Z} T^{np}(0,p)_{I_{\halb}}$. Further,
(\ref{eq:lastret}) implies that $\gamma$ does never intersect itself.

$\gamma$ can be strobed at $\theta = 0$ by setting $x_n := \gamma(n) \ 
\forall n \in \Z$. Let $\Lambda := \{ \xn \mid n \in \Z \}$. Then the
map $ \Lambda \ra \Lambda, \ \xn \mapsto x_{n+1}$ is order-preserving
and bijective, but this means that the sequence \nfolge{\xn} is
combinatorically equivalent to the orbit of a circle homeomorphism.
In addition, to each $n \in \Z$ there exists a unique
number $a(n) \in \Z$ which satisfies
\[
    (0,\xn) \in T^{a(n)p}(0,p)_{J} \ \textrm{ and } \ a(n) \in A := \{
    a\in \Z \mid |a\omtil \mod 1| \leq \frac{\omtil}{2} \}.  
\] 
$(0,\xn)$ may be contained in several successive iterates of
$(0,p)_{I_{\halb}}$, but the second requirement makes the choice of
$a(n)$ unique (note that $|a\omtil \bmod 1| = \halb \omtil$ is not
possible since $\omtil$ is irrational). If we keep in mind the way in
which the iterates of $(0,p)_{I_{\halb}}$ are attached to and moving
along $\gamma$, it is obvious that the map $a : \Z \ra A$ defined in
this way is bijective, monotone and symmetric (i.e.\ $a(-n) = -a(n)$).

\ \\
Similar as before we will call $n \in \Z$ a \textit{closest
  intersection time} (in order to distinguish it from the closest
return times), whenever \ $x_{-n} < \xnull < \xn$ and $\nexists k \in
\Z \setminus\{0\}: |k| < |n|, \ x_{-n} < x_k < \xn$ \ or \ $\xn <
\xnull < x_{-n}$ and $\nexists k \in \Z\setminus \{0\}: |k| < |n|, \ 
x_n < x_k < x_{-n}$. Again we want to apply Lemma \ref{lem:closestret}
to conclude that there can only be finitely many such $n$.

To that end, consider the rectangle $W' = I' \times K$, where $I'$
is the symmetric middle part of $I$ with lenght exactly $\frac{3}{2}
\omtil$. Note that $I' \subseteq I_{\halb}$, because $|\omtil| =
|p\omega \bmod 1| \leq \viertel |I|$ as $p \in N(\halb)$. $W'$ is a
wandering rectangle for $T^p$, and we now want to look at the closest
return times of $W'$ with respect to $I'_{\drittel}$. However, we have 
already described them quite neatly: First of all, the set of return
times $N'(\drittel)$ is exactly the set $A$ defined above. As
$T^{a(n)}W'$ is contained in $T^{a(n)}(0,p)_{I_{\halb}}$, the ordering 
of these sets over $I'(\drittel)$ coincide and are determined by the
ordering of the points $\xn \in \Lambda$. Thus $k \in \Z$ is a closest 
return time (with respect to $W'$ and $I'_{\drittel}$) if and only if
$k = a(n)$ for some closest intersection time $n$. As there cannot be
infinitely many closest return times by Lemma \ref{lem:closestret},
the same is true for the closest intersection times. 

Let $q$ be the maximum of the closest intersection times. First
assume $q = 1$ and (w.l.o.g.) $x_{-1} < \xnull < x_1$. Then
\nfolge{\xn} is strictly monotonically increasing in the sense that 
\[
    \xnull < \xeins < \ldots < \xn < x_{-n} < \ldots x_{-1} < \xnull \ 
    \ \forall n \in \N \ .
\]
A sequence of graphs which increases monotonically in the same sense
can be obtained by 
\[
    \varphi_n(\theta) := \gamma(\thhat) \ \textrm{ with }
    \pi(\thhat) = \theta \ \textrm{  and } \thhat \in [n,n+1) 
\]
(such that $\varphi_n$ can be identified with $\gamma_{|[n,n+1)}$ in a
very natural way). We have
\[
    \Tth^p(\varphi_n(\theta)) \ = \ \left\{ \begin{array}{ll}
    \varphi_n(\theta + \omtil) & \textrm{if } \theta \in [0,1-\omtil)
    \\
    \varphi_{n +1}(\theta + \omtil)  &\textrm{if } \theta \in [1-\omtil,0)
    \end{array} \right.
\]
Further the graphs $\varphi_n$ are continuous, except at $\theta = 0$ where
they are continuous to the right and $\lim_{\theta \nearrow 0}
\varphi_{n}(\theta) = \varphi_{n+1}(0)$. From all this it follows
easily that $U := \bigcup_{n\in\Z} (\varphi_{-n},\varphi_{n})$ is an
$1,1$-invariant open tube for $T^p$ (and thus some $p',1$-invariant
tube for $T$ with $p' \leq p$ must exist as well). 

Now if $q>1$ we can lift both $\gamma$ and $T^p$ to the ``blown-up''
torus $\R / q\Z \times \kreis$. If we then strobe the lift
$\bar{\gamma}$ at $\theta = 0$, the set $\bar{\Lambda}$ we obtain
contains exactly the points $x_{nq} \ (n\in\Z)$ and the construction
works as before. The $p',1$-invariant open tube $\bar{U}$ on $\R/q\Z$
will then project down to a $p',q$-invariant open tube for the
original system. That this projected tube will not ``intersect
itself'' (i.e.\ the projection restricted to the tube is injective)
follows from the fact that for no $k,n\in\Z$ $x_k$ can be contained in
$(x_{nq},x_{(n+1)q})$. Thus 
$\bar{U}_{\theta=0} = \bigcup_{n\in\Z} (x_{-nq},x_{nq})$ and
$\bar{U}_{\theta=j} = \bigcup_{n\in\Z} (x_{-(nq+j)},x_{nq+j})$ are
disjoint $\forall j = 1 \ld q-1$. As strobing $\gamhat$ at $\theta
= 0$ was arbitrary ($\bar{U}$ could have been obtained in the same way
by strobing at any other $\theta$), the same is true for any $\theta
\in \R/q\Z$, i.e.\ $\bar{U}_\theta$ and $\bar{U}_{\theta+j}$ are
disjoint $\forall j = 1 \ld q-1$.

\ \\
\underline{Step 2:} \ Absence of regular invariant graphs implies
transitivity. \\
We have to show that for all open sets $U,V \subseteq \ntorus \
\exists n\in\N: T^{-n}U \cap V \neq \emptyset$. Of course it suffices
to restrict to the case where $U$ and $V$ are rectangles (in
particular connected). First suppose
\begin{equation}
  \label{eq:disjoint}
  T^nU \cap V = \emptyset \ \forall n \in \Z \ .
\end{equation}
As $U$ is non-wandering we have\ $T^mU \cap U \neq \emptyset$ for some
$m\in \N$. Thus $\tilde{U} := \ncup T^{nm}U$ will be a kind of
``infinite open tube'', similar to the situation in Step 1. But in
this case, as $U$ is non-wandering for any iterate of $T$, the tube
has to intersect itself after winding around the circle a certain
number, say $q$, of times. By passing to the $q$-fold torus again if
necessary, we can assume $q=1$. This means that for any point $\thx
\in \tilde{U}$ there exists a simple closed curve $\gamma_{\thx}$
passing through \thx, which is completely contained in $\tilde{U}$.

Now $\tilde{U}$ is a $T^m$-invariant set, and by using some multiple
$p$ of $m$ we can repeat this construction to obtain a $T^p$-invariant 
set $\tilde{V} := \ncup T^{np}V$ with exactly the same
properties. (\ref{eq:disjoint}) implies $\tilde{U} \cap \tilde{V} =
\emptyset$. 

$\tilde{U}$ and $\tilde{V}$ are already very close to $p,q$-invariant
open tubes, but they may still have some ``holes''. However, points in 
$\tilde{U}$ and $\tilde{V}$ cannot alternate on the fibres: Suppose
$u_1,u_2 \in U_\theta$ and $v_1,v_2 \in V_\theta$. The curves
$\gamma_{(\theta,u_1)}$ and $\gamma_{(\theta,u_2)}$ divide the torus
into at least two connected components. As $\tilde{V}$ is connected
(which is obvious from the construction) it must be completely
contained in one of them, i.e.\ $v_1$ and $v_2$ must be contained in
the same interval $(u_1,u_2)$ or $(u_2,u_1)$. Therefore $u_1 < v_1 <
u_2 < v_2 < u_1$ is not possible. With this it is easy to see, that
\[
    \bar{U} := \{ \thx \mid \exists u_1,u_2 \in \tilde{U}_\theta : x \in
               [u_1,u_2] \textrm{ and } [u_1,u_2] \cap V_\theta =
               \emptyset \} 
\]
is a $1,1$-invariant open tube with respect to $T^p$, contradicting
the assumption that there are no regular invariant graphs.

\ \\
Thus we have $\exists n \in \Z: T^nU \cap V \neq \emptyset$. If $n<0$
we are done. Otherwise $U \cap T^{-n}V$ is non-wandering with respect
to $T^n$, i.e.\ $\emptyset \neq (U \cap T^{-n}V) \cap T^{kn}(U \cap
T^{-n}V) \subseteq U \cap T^{(k-1)n}V$ for some $k \geq 1$. Therefore
$T^{-(k-1)n}U \cap V \neq \emptyset$.
   
\qed

\ \\
Recall that for a single circle diffeomorphism transitivity implies
minimality. For quasiperiodically driven circle diffeomorphisms we do
not expect this to be true. But the preceding proof can be modified a
bit to give information about the structure of all minimal sets of the
system. Indeed, what we proved is that $\ncup T^{-n}U$ is dense in
$\ntorus$ for each open set $U\subseteq\ntorus$. Equivalently, the set
of points $\thx$ whose orbit never enters $U$ is closed and nowhere
dense so that the set $N$ of points $\thx$ whose orbit is not dense is
meager (\emph{i.e.}  of first Baire category). If $T$ is not minimal,
then obviously $N$ contains every minimal subset of $\ntorus$. The structure of such minimal sets is described in the following theorem.

\begin{thm}\label{thm:disconnectedness}
  Suppose that $T\in\Tbv$ has no regular invariant graph (so that it
  is transitive) but that it is not minimal, and let $M$ be a minimal
  subset of $\ntorus$. Then each connected component of $M$ is
  contained in a single fiber $\pi_1^{-1}(\theta)$.
\end{thm}
\proof\\
Fix some rectangle $U$ in the complement of a minimal set $M$. As in
Step~2 of the preceding proof the ``infinite open tube'' $\tilde
U=\ncup T^{nm}U$ intersects itself after winding around the torus $q$
times where we can assume as above that $q=1$, and there exists a
simple closed curve $\gamma$ winding around the torus which is
completely contained in $\tilde U$ so that it has distance $\delta_M>0$ to
the minimal set $M$. If we construct this curve as in Step~1 of the
preceding proof, except that in the last step of winding around once
we close it, then $\Gamma\cap T^m\Gamma\neq\emptyset$.  Denote by $Q$
the complement of $\Gamma\cup T^m\Gamma$ in $\ntorus$ and observe that
$M\subseteq Q$.  

Now suppose for a contradiction that $C$ is a connected component of
the minimal set $M$ such that $\pi_1(C)$ is a nontrivial subinterval
of $\kreis$. Denote its length by $\delta_C$. Let $Q'$ be any
connected component of $Q$ and assume that $Q'\cap M\neq\emptyset$.
Then $Q'\cap M$ is closed and so is $\pi_1(Q'\cap M)$. We claim that
$\pi_1(Q'\cap M)$ is also open: Let $\thx\in Q'\cap M$ and fix some
point $(\eta,y)\in C$ such that $\eta$ is in the middle part of length
$\delta_C/2$ of the interval $\pi_1(C)$. As the orbit of
$(\eta,y)$ is dense in $M$ there is some iterate $k$ such that
$d(T^k(\eta,y),\thx)<\min\{\delta_M,\delta_C/4\}$ and this implies at
once that $T^kC\subseteq Q'\cap M$ and that $\pi_1(Q'\cap
M)\subseteq\kreis$ contains an interval neighbourhood of $\theta$.
Hence the nonempty set $\pi_1(Q'\cap M)$ is closed and open so that
indeed $\pi_1(Q'\cap M)=\kreis$. But since $\Gamma\cap
T^m\Gamma\neq\emptyset$, the complement $Q$ of $\Gamma\cup T^m\Gamma$
can have at most one such connected component.

Equipped with this information we will now construct an invariant open
tube which contradicts the assumption that $T$ has no regular
invariant graph: Define $\hat\phi,\hat\psi:\kreis\to\R$ by
\begin{displaymath}
  \begin{split}
  \hat\phi(\theta)
  &:=
  \inf\{t>0:\gamma(\theta)-t\in M_\theta\}\\
  \hat\psi(\theta)
  &:=
  \inf\{t>0:\gamma(\theta)+t\in M_\theta\}    
  \end{split}
\end{displaymath}
As $\gamma$ is a continuous curve and as $M$ is closed, the function
$\hat\phi$ is u.s.c. and $\hat\psi$ is l.s.c. and $\gamma-\phi$ and
$\gamma+\psi$ are u.s.c. respectively l.s.c. continuous functions from
$\kreis$ to $\kreis$. $(\gamma-\phi,\gamma+\psi)$ is an open tube
disjoint from $M$, and we will show that it is $T^m$-invariant: By
definition of $\hat\phi$,
\begin{displaymath}
  \begin{split}
  \Tth^m(\gamma(\theta)-\phi(\theta))
  &=
  \Tth^m(\gamma(\theta))-
  \inf\{t>0:\Tth^m(\gamma(\theta))-t\in M_{\theta+m\omega}\} \quad \text{whereas}\\    
  \gamma(\theta+m\omega)-\phi(\theta+m\omega)
  &=
  \gamma(\theta+m\omega)-\inf\{t>0:\gamma(\theta+m\omega)\in M_{\theta+m\omega}\}\ .
  \end{split}
\end{displaymath}
But the two right hand sides have the same value, because $M$ is
contained in a single connected component of the complement of
$\Gamma\cup T^m\Gamma$ so that for any two points $u,v\in
M_{\theta+m\omega}$ it is impossible that
$\gamma(\theta+m\omega)<u<\Tth^m(\gamma(\theta))<v<\gamma(\theta+m\omega)$.

\qed

\ \\
The kind of total disconnectedness of minimal sets in
$\theta$-direction that is expressed in
Theorem~\ref{thm:disconnectedness} is in some contrast to the
following simple observation.

\begin{lem}
  Let $T\in\Tbv$ and let $\emptyset\neq M$ be a minimal strict subset
  subset of $\ntorus$. Then every point $\thx\in M$ has the following
  property: If $U$ is any neighbourhood of $\thx$, then $\pi_1(U)$
  contains a nontrivial interval which contains $\theta$. (It may
  happen that $\theta$ is an endpoint of this interval.)
\end{lem}
\proof \\
Denote by $G$ the set of theose $\thx\in\ntorus$ which have a
neighbourhood $U$ such that $\pi_1(U\cap M)$ is nowhere dense in
$\kreis$.  This set is obviously open and invariant under $T$. Hence
$M\setminus G$ it is either empty or all of $M$. Suppose first that
$M\subseteq G$.  Then the compact set $M$ can be covered by finitely
many sets $U$ for which $\pi_1(U\cap M)$ is nowhere dense, and it
follows that $\pi_1(M)$ is a nowhere dense closed subset of $\kreis$
invariant under the irrational rotation by $\omega$. Hence
$\pi_1(M)=\emptyset$ in contradiction to the assumption
$M\neq\emptyset$. Therefore $M\cap G=\emptyset$ so that each point
$\thx\in M$ has the property claimed in the lemma.

\begin{figure}[H!]
\vspace*{-6eM}
  \centering
  \epsfig{file=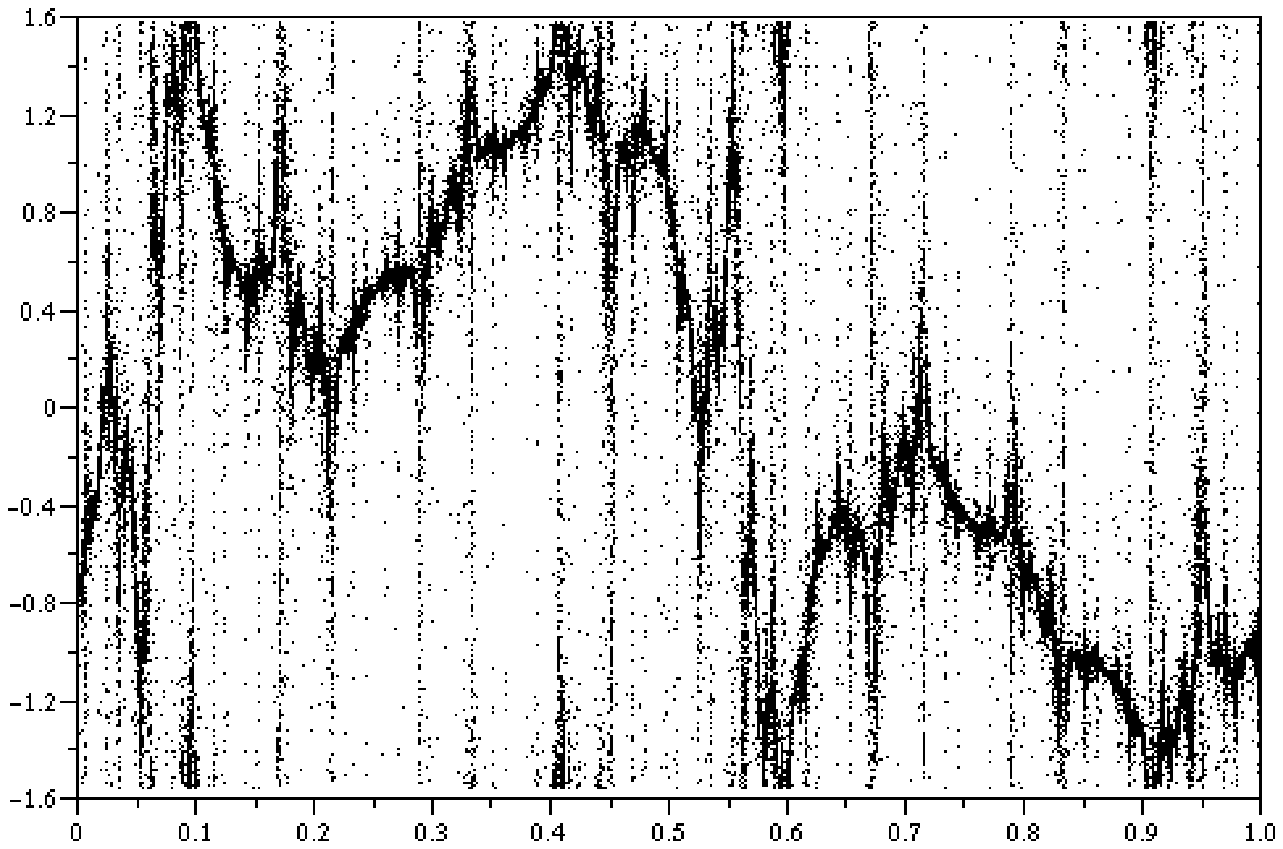,height=10cm,width=16cm}
  \caption{A trajectory of $T$}
  \vspace*{-7.5cm}
  \begin{sideways}
    \hspace*{5cm}$\arctan(x)$
  \end{sideways}
  \quad\hspace*{12cm}$\theta$
  \label{fig:1}
\end{figure}

\begin{figure}[H]
  \epsfig{file=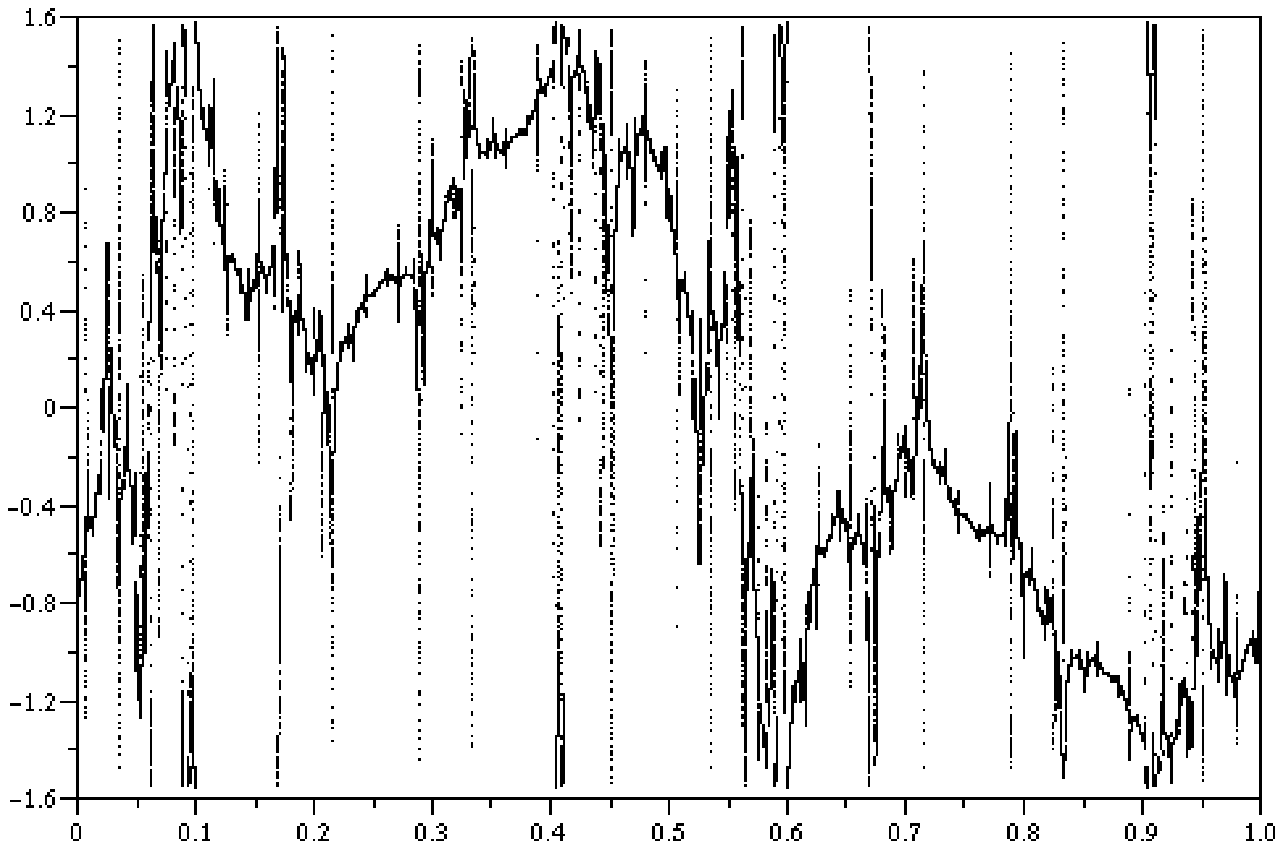,height=10cm,width=16cm}
  \caption{The invariant graph}
  \vspace*{-7.5cm}\hspace*{0.6cm}
  \begin{sideways}
    \hspace*{5cm}$\arctan(x)$
  \end{sideways}
  \quad\hspace*{12cm}$\theta$
  \label{fig:2}
\end{figure}

\qed

\ \\
We close this section with an example which seems a good candidate for
transitive but nonminimal behaviour. It is the \emph{critical Harper
  map} 
\begin{displaymath}
  T\thx=\left(\theta+\omega,\frac{-1}{x-E+\lambda\cos(2\pi\theta)}\right)
  \quad\text{with $E=0$ and $\lambda=2$}
\end{displaymath}
with fiber maps $\Tth$ acting on $P^1\R$. After a change of
co-ordinate $x'=\arctan(x)$ the map $T$ belongs to the class $\Tbv$.
The map $T$ has Lyapunov exponent zero
\cite{bourgain/jitomirskaya:2002} and rotation number $\frac{1}{2}$.
\footnote{This follows easily from the symmetry
  $T_{\theta+\frac{1}{2}}(x)=-\Tth(-x)$ of the map on $\kreis \times
  \bar{\R}$:  This implies that whenever $\That : \kreis \times \R \ra
  \kreis \times \R$ is a lift of $T$, then so is the map
  $\overline{T}$ defined by $\overline{T}_\theta(x) :=
  -\That_{\theta+\halb}(-x)$. Thus the rotation number must be $0$ or
  $\halb$, and this is true for any $\lambda$. Now $\rho(T)=\halb$ for
  $\lambda=0$ and the rotation number depends continuously on
  $\lambda$.}  Numerical evidence suggests that $T^2$ has no invariant
tube, and invariant tubes for higher powers $T^q$ which are not
invariant for any smaller power of $T$ are uncompatible with rotation
number $\frac{1}{2}$.  Therefore $T$ should be transitive. Numerical
evidence again, but also the general classification of ergodically
driven M\"obius maps \cite{thieullen:1997,arnold/cong/oseledets:1999}
together with the close relation of the dynamics of this family of
maps with spectral properties of the almost Mathieu operator
\footnote{See
  \cite{ketoja/satija:1997,prasad/ramaswamy/satija/shah:1999} for this
  relation and \cite{jitomirskaya:1999} for the relevant spectral
  properties.}  suggests that $T$ has a unique invariant probability
measure $\mu(A)=\int_{\kreis} {\bf 1}_A(\theta,\phi(\theta))\,d\theta$
where $\phi:\kreis\to\kreis$ is a measurable invariant graph.  (This
is called the parabolic case in \cite{thieullen:1997}.) Then the
topological support of $\mu$ would be the only minimal set.
Figure~\ref{fig:1} shows the plot of a trajectory of length $10^5$,
and Figure~\ref{fig:2} displays the result of a numerical
reconstruction of the graph $\phi$ based on the assumption that the
map is indeed parabolic. More details can be found in the forthcoming
note \cite{datta/jaeger/keller/ramaswamy:2003}.

%  \vspace*{-10eM}
\newpage
%
%
%                Literaturverzeichnis
%
%

\bibliography{snaphysics,qpfs}  \bibliographystyle{alpha}

\end{document}